\UseRawInputEncoding
\documentclass[12pt]{amsart}
\usepackage{amsmath,
			amssymb,
			amsthm,
			amsfonts,
			mathrsfs, 
			enumerate, 
			mathtools} 
\usepackage[all]{xy}
\usepackage[normalem]{ulem}
\usepackage{url}
\usepackage{color}

\usepackage{pifont}

\usepackage{mathabx}

\usepackage{hyperref}
\usepackage{graphicx}

\usepackage{txfonts}

\newtheorem{theorem}{Theorem}

\newtheorem{proposition}{Proposition}[section]
\newtheorem{lemma}[proposition]{Lemma}
\newtheorem{corollaire}[proposition]{Corollary}

\newtheorem{definition}[proposition]{Definition}
\newtheorem{example}[proposition]{Example}

\theoremstyle{definition}
\newtheorem{remarque}[proposition]{Remark}

\newcommand{\bl}{\begin{lemma}}
\newcommand{\bp}{\begin{proposition}}
\newcommand{\bt}{\begin{theorem}}
\newcommand{\bc}{\begin{corollaire}}
\newcommand{\be}{\begin{equation}}
\newcommand{\bee}{\begin{equation*}}
\newcommand{\bd}{\begin{definition}}
\newcommand{\bdp}{\begin{definitionproposition}}
\newcommand{\bex}{\begin{example}}
\newcommand{\br}{\begin{remarque}}
\newcommand{\bpr}{\begin{proof}}

\newcommand{\el}{\end{lemma}}
\newcommand{\ep}{\end{proposition}}
\newcommand{\et}{\end{theorem}}
\newcommand{\ec}{\end{corollaire}}
\newcommand{\ee}{\end{equation}}
\newcommand{\eee}{\end{equation*}}
\newcommand{\ed}{\end{definition}}
\newcommand{\edp}{\end{definitionproposition}}
\newcommand{\eex}{\end{example}}
\newcommand{\er}{\end{remarque}}
\newcommand{\epr}{\end{proof}}

\newcommand\Z{\mathbb Z}

\newcommand\nada{\phantom .}

\renewcommand\1{\hbox{\ding{192}}}
\renewcommand\2{\hbox{\ding{193}}}
\newcommand\3{\hbox{\ding{194}}}
\newcommand\4{\hbox{\ding{195}}}

\newcommand{\sbat}{S^1}

\newcommand{\stres}{S^3}
\newcommand{\zdos}{{\Z}_{_2}}

\newcommand{\lau}[4]{{#1}^{^{#2}}_{_{#3}}{\left( #4 \right)}}
\newcommand{\coho}[3]{{#1}^{^{#2}}{\left( #3 \right)}}
\newcommand{\homo}[3]{{#1}_{_{#2}}{\left( #3 \right)}}

\newcommand
{\SucExaRef}[5]
{
\begin{equation}\label{#1}
\xymatrix@C=0,4cm
{
\cdots \ar[r] & \ #2\ar[r] &#3\ar[r] &  #4\ar[r] & #5 \ar[r] & \cdots ,
}
\end{equation}
}

\newcommand
{\SucExa}[4]
{
\begin{equation*}
\xymatrix@C=0,4cm
{
\cdots \ar[r]  &  #1\ar[r]  &  #2\ar[r]  & #3\ar[r] & #4 \ar[r] & \cdots ,
}
\end{equation*}
}

\def\per{\overline}

\newcommand{\menos}{\backslash}
\newcommand{\rondp}{\raise1pt\hbox{\tiny $\circ$}}
\newcommand\phii{{\raise2pt\hbox{$\varphi$}}}

\newcommand\Om{\Omega}  
\newcommand\om{\omega}

\newcommand\chii{\raise2pt\hbox{$\chi$}}

\title{Smith-Gysin sequence}

\date{\today}

\author{Jos\'{e} Ignacio Royo Prieto
}
\address{Matematika Saila\\ Zientzia eta Teknologia Fakultatea\\ University of the Basque Country UPV/EHU\\ Barrio Sarriena s/n\\ 48940 Leioa\\Spain.
\thanks{Partially supported by Ministerio de Ciencia, Spain, grant PID2019-105621GB-I00.}
}

\email{joseignacio.royo@ehu.eus}

\author{Martintxo Saralegi-Aranguren}
\address{Laboratoire de Math{\'e}matiques de Lens\\  
      EA 2462 \\
      Universit\'e d'Artois\\
         SP18, rue Jean Souvraz\\
          62307 Lens Cedex\\
         France}
\email{martin.saraleguiaranguren@univ-artois.fr}

\author{  
Robert Wolak}
\address{Instytut Matematyki. Uniwersytet Jagiellonski. \\
Stanislawa Lojasiewicza 6, 30\\
Krakow, Poland \thanks{Partially supported by the KBN grant 2 PO3A  021 25.}}
\email{robert.wolak@im.uj.edu.pl.}  

\keywords{Long exact sequence, $S^3$-actions.}
\subjclass[2010]{Primary 57S15; Secondary 55N10	.}
\thanks{
The authors acknowledge that the research cooperation was funded by the program Excellence Initiative Ð Research University at the Jagiellonian University in Krakow within the framework of the research group Reeb-Reinhart 2022.}

\begin{document}

\begin{abstract}
Starting with a manifold $M$ and a semi-free action of $S^3$ on it, we have the Smith-Gysin sequence:
$$
\cdots \to  H^{*}( M)
\to
  H^{*-3}(M/S^3, M^{S^3})  \oplus  H^{*} (M^{S^3})
\to
H^{*+1}(M/S^3, M^{S^3})
\to
 H^{*+1}(M)
\to
\cdots
$$
In this paper, we construct a Smith-Gysin sequence that does not require the semi-free condition. This sequence includes a new term, referred to as the "exotic term," which depends on the subset $M^{\sbat}$:
$$
\cdots \to 
H^{*}(M)
\to H^{*-3} (M/S^3, \Sigma/S^3)  \oplus  H^{*}(M^{S^3}) \oplus \left( H^{*-2}(M^{\sbat})\right)^{-\zdos}
\to
H^{*+1}(M/S^3,M^{S^3})
\to
H^{*+1}(M)
\to \cdots
$$
 Here, $\Sigma \subset M$ is the subset of points in $M$ whose isotropy groups are infinite. The group $\mathbb{Z}_2$ acts on $M^{S^1}$ by $j \in S^3$.
 \end{abstract}

\maketitle

Let us consider a smooth action of $S^3$ on a manifold $M$. The cohomological relationship between the different actors involved in this action can be expressed through a long exact sequence known as the Gysin sequence. When the action is free, we have the following exact sequence:
\SucExa{\coho H {*} M} {\coho H {*-3} {M/S^3}} {\coho H {*+1} {M/S^3} } {\coho H {*+1} M }
(see for example \cite{MR0336651,MR0413144}). Here, $M/S^3$ denotes the orbit space which is a smooth manifold in this context.
Suppose that the action is free outside the fixed point set $M^{S^3}$ (\emph{semi-free action}).
When the action is semi-free, we get 
\SucExaRef{gys3SemFr}{\coho H {*} M} {\coho H {*-3} {M/S^3, M^{S^3}} } {\coho H {*+1} {M/S^3} } {\coho H {*+1} M }
(see for exemple \cite[Corollary]{MR1184085}). Finally, in the general case, an exotic term appears:
\SucExa
{ \coho H {*} M }{ \coho H {*-3} {M/S^3, \Sigma/S^3} \oplus \left( \coho H {*-2} {M^{\sbat}}\right)^{-\zdos}}{ \coho H {*+1} {M/S^3}}{ \coho H {*+1} M}
(cf. \cite{MR3119667}).
 Here, $\Sigma \subset M$ is the subset of points in $M$ whose isotropy groups are infinite.  The group $\zdos$ acts on $M^{\sbat}$ by $j \in \stres$.

On the other hand, in the case of a semifree action, we also have the Smith-Gysin sequence
\SucExaRef{gys3ClasSG}
{
\coho H {*} M}{ \coho H {*-3} {M/S^3, M^{S^3}} \oplus \coho H {*} {M^{S^3}}}{
\coho H {*+1} {M/S^3, M^{S^3}}}{ \coho H {*+1} M }
(see, for example, \cite[Theorem 2.2]{Smith-Gysin} and \cite[Exercise 12 on page 169]{MR0413144}).

In this paper, we extend the scope of the Smith-Gysin sequence to any smooth action of $S^3$ on a manifold $M$, obtaining the long exact sequence
\SucExa
{
\coho H {*} M 
}{
 \coho H {*-3} {M/S^3, \Sigma/S^3} \oplus \left( \coho H {*-2} {M^{\sbat}}\right)^{-\zdos} \oplus  \coho H *{M^{S^3}} 
}{
 \coho H {*+1} {M/S^3,M^{S^3}} 
 }{
  \coho H {*+1} M
  }

In what follows, we consider a smooth action $\Phi \colon S^3 \times M \to M$, where $M$ is a second-countable, Hausdorff, smooth manifold of dimension $m$ and without boundary. We further assume that the action is effective. For definitions and properties related to actions of compact Lie groups, we refer the reader to \cite{MR0413144}.
The cohomology $\coho H* -$ refers to  the singular cohomology withe real coefficients.

\section{Verona's differential forms}

The action $\Phi$ induces the following invariant filtration\footnote{We refer to \cite{MR3119667} for the notions appearing in this presentation.} on $M$:
$
F  \subset \Sigma \subset M,
$
where 

\nada \hfill  $F = \{x \in M \mid \dim S_x^3  =3\}$ \hfill $\Sigma = \{x \in M \mid \dim S_x^3  = 1,3\}$ \hfill \nada

\noindent Notice that $F$ is the fixed point set $M^{S^3}$.

 The complex of controlled forms $\lau \Om *VM$ (or Verona's forms), which is composed of differential forms defined on $M\setminus\Sigma$, computes the usual cohomology $\mathrm{H}^{*}(M)$. In fact, a Verona's form can be represented as a triple $(\omega, \omega_0, \omega_1)$ where

\smallskip

\nada \hfill $\omega \in \coho \Omega * {M\menos \Sigma}$ \hfill $\om_0\in \coho \Om * {\Sigma\menos F}$
 \hfill  $\om_1\in \coho\Om * F$ \hfill \nada

\smallskip

\noindent such  that  $\omega = \omega_0 = \omega_1$ near $F$, and
$\omega = \omega_0$ near $\Sigma \menos F$ (see \cite[Definition 1.2]{MR3119667} for the precise statement).
In this work, we use the relative complex
$$
\lau \Om *V{M,F} = \{ \omega \in \lau \Om * V M \mid \omega_1 =0\}.
$$
\begin{proposition}\label{Pbat}
Let $\Phi \colon \stres \times M \to M$ be a smooth action. The long sequence
\begin{equation}\label{bat}
\xymatrix{
0\ar[r] &  \lau \Om * V{M,F} \ar@{^(->}[r] & \lau \Om * V M \ar[r]^{\rho_F} & \coho \Om * F \ar[r] \ar[r] & 0,
}
\end{equation}
where $\rho_F(\omega) = \omega_1$,
is exact. Moreover, we have $\coho H* {\lau \Om \cdot V{M,F}} = \coho H * {M,F}$.
\end{proposition}
\bpr
The long sequence is exact if we prove that the map $\rho_F$ is a surjective map. It suffices to proceed as in \cite[Lemma 1.2]{MR3119667}. On the other hand, Verona's forms are a particular case of intersection differential forms. In fact, we have $\lau \Om *{V}{M,F} =\lau \Om*{\per p}M$, where $\per p$ is the perversity taking the value $-\infty$ on the strata of $F$ and $0$ on the strata of $\Sigma\setminus F$. The result comes from \cite[Proposition 3.2.3]{MR2210257}. 
\epr

The orbit space $M/\stres$ is a stratified pseudomanifold and its cohomology can be also computed by using differential forms. More specifically, the complex 
$$
\lau \Om * V {M/\stres} = \{ \omega \in \lau \Om * V{M} \mid i_X\omega = i_X d \omega = 0\ \ \ \forall X \in \homo{\mathfrak{X}}{\Phi}{M}\}
$$
computes the cohomology $\coho H * {M/\stres}$ (cf. \cite{MR290375}). 
Here, we denote by $\homo{\mathfrak{X}}{\Phi}{M}$
the family of  vector fields of $M$ tangent to the orbits of $\Phi$.
We also use the relative complex
$$
\lau \Om * V {M/\stres,F} = \{ \omega \in \lau \Om * V{M,F} \mid i_X\omega = i_X d \omega = 0\ \ \ \forall X \in \homo{\mathfrak{X}}{\Phi}{M}\}
$$

Proceeding as in the previous proposition, we obtain the following result:
\begin{proposition}\label{Pbi}
Let $\Phi \colon \stres \times M \to M$ be a smooth action. The long sequence
\begin{equation}\label{bi}
\xymatrix{
0\ar[r] &  \lau \Om * V{M/\stres ,F} \ar@{^(->}[r] & \lau \Om * V {M/\stres}  \ar[r]^-{\rho'_F} & \coho \Om * F \ar[r] \ar[r] & 0,
}
\end{equation}
where $\rho'_F(\omega) = \omega_1$,
is exact. Moreover, we have $\coho H* {\lau \Om \cdot V{M/\stres,F}} = \coho H * {M/\stres,F}$.
\end{proposition}

\section{Smith-Gysin sequence}

In order to obtain the Smith-Gysin sequence, we use the notion of a braid.

 \begin{definition}
 Let us consider six chain complexes $A^*,B^*,C^*,D^*,E^*$ and $F^*$.
 A {\em braid} is a diagram of chain maps of the form
 
 \bee
\scalebox{1}
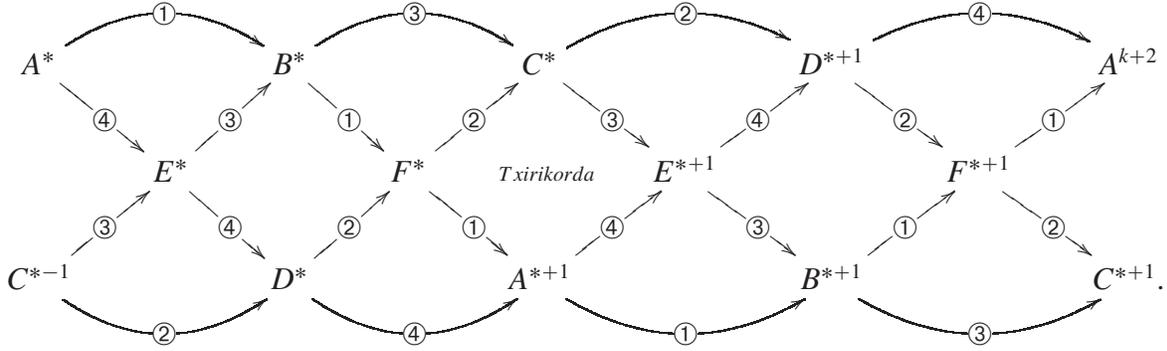
{
 \xymatrix{
A^*  \ar[rd] |\4\ar@/^{7mm}/[rr] |\1 &
&
B^*\ar[rd] |\1\ar@/^{7mm}/[rr] |\3&
&
 C^* \ar[rd] |\3 \ar@/^{7mm}/[rr] |\2&
&
D^{*+1} \ar[rd] |\2 \ar@/^{7mm}/[rr] |\4&
&
A^{k+2} \\
&
E^* \ar[ur]  |\3\ar[dr] |\4 &
&
 F^*\ar[ur]  |\2 \ar[dr] |\1 \ar@{} [rr]| {Txirikorda}&
&
E^{*+1} \ar[ur] |\4\ar[dr] |\3 &
&
F^{*+1} \ar[ur] |\1\ar[dr] |\2&
\\
C^{*-1}   \ar[ur] |\3\ar@/_{7mm}/[rr] |\2&
&
D^* \ar[ur] |\2\ar@/_{7mm}/[rr]  |\4&
&
A^{*+1} \ar[ur] |\4\ar@/_{7mm}/[rr] |\1&
&
B^{*+1} \ar[ur] |\1\ar@/_{7mm}/[rr] |\3&
&
C^{*+1}.
}}
\eee

\smallskip

 \noindent It is a {\em commutative braid} when all the triangles and diamonds are commutative. If the long sequences $\1$, $\2$, $\3$ and $\4$ are exact we say that braid is an \em{exact braid}.
 \end{definition}

An exact and commutative braid has the following property: the long sequence shown below is exact (see, for example, \cite[pag. 39-41]{MR0178036}).
\be\label{B1}
\xymatrix@C=1cm{
\cdots \ar[r] &
E^* \ar[r]^-{( \3 , \4)} &
B^* \oplus D^* \ar[r]^-{\1 - \2}& F^* \ar[r]^-{\3 \2} & E^{*+1} \ar[r] & \cdots
}
\ee

The main result of this work is the following.

\begin{proposition}
Consider a smooth action $\Phi \colon S^3 \times M \to M$, and let $\Sigma = \{ x \in M \mid \dim S^3_x=1,2\}$. The Smith-Gysin sequence associated to this action is the following long exact sequence:
$$
\xymatrix@C=0,33cm
{
\cdots \ar[r] &
\coho H {*-1} M 
\ar[r]  &\coho H {*-4} {M/S^3, \Sigma/S^3}  \oplus \coho H {*-1} {M^{S^3}} \oplus \left( \coho H {*-3} {M^{\sbat}}\right)^{-\zdos}
\ar[r]  &
 \coho H {*} {M/S^3,M^{S^3}} 
\ar[r]  &
 \coho H {*} M
\ar[r]  &\cdots}
$$
Here, the group $\mathbb{Z}_2$ acts on $M^{\sbat}$ by $j \in \stres$.
We also have the long exact sequence
$$
\xymatrix@C=0,35cm
{
\cdots \ar[r] &
 \coho H {*}{ M/\stres, M^{\stres}}  \ar[r]  &   
\coho H {*} {M/{S^3}}  \oplus  \coho H {*} { M, M^{\stres}}  \ar[r]  &
 \coho H {*} {M} \ar[r]  &
 \coho H {*+1}{ M/\stres, M^{\stres}} \ar[r]  & 
 \cdots,
}
$$
\end{proposition}
\bpr
The long exact sequences \eqref{bat} and \eqref{bi} can be arranged in the following commutative diagram: 
$$
\scalebox{.9}{
\xymatrix@=.5cm{
  & 0 \ar[d]& 0 \ar[d] &  & &
  \\
    0 \ar[r]& \lau \Om * V {M/\stres,{M^{\stres}} } \ar@{^(->}[r] \ar@{^(->}[d] & \lau \Om * V {M/\stres}  \ar@{^(->}[d]\ar[r]&  
\coho \Om * {M^{\stres}} \ar[r]  \ar@{=}[d] &  0 & \3
  \\
    0 \ar[r]& \lau \Om * V {M,{M^{\stres}} } \ar[r] \ar@{^(->}[d] & \lau \Om * V M  \ar[d]  \ar[r]&  
\coho \Om * {M^{\stres}} \ar[r] &  0 & \2
    \\  
   & \lau Q* \Phi M\ar@{=}[r] \ar[d]&\lau Q* \Phi M\ar[d]& & &\\ 
  & 0&0 && & \\
   & \4&\1 && &  
}}
$$

\bigskip

The cohomology  $\coho H * {M,M^{\stres}}$ have been computed in \cite[Lemma 2.1]{MR3119667}. The result is given by
$$
\coho H {*-3} {M/S^3, \Sigma/S^3}  \oplus \left( \coho H {*-2} {M^{\sbat}}\right)^{-\zdos}.
$$
This diagram leads to the following braid  (cf. Propositions \ref{Pbat} and \ref{Pbi}).

 \bee
\scalebox{.55}{
 \xymatrix{
&
\coho H * {M/\stres} \ar[rd] |\1\ar@/^{7mm}/[rr] |\3&
&
\coho H * {{M^{\stres}} }  \ar[rd] |\3 \ar@/^{7mm}/[rr] |\2&
&
\coho H {*+1} {M,{M^{\stres}} }  \ar[rd] |\2 \ar@/^{7mm}/[rr] |\4&
&
\coho H {*-2} {M/S^3, \Sigma/S^3}  \oplus \left( \coho H {*-1} {M^{\sbat}}\right)^{-\zdos} 
\\
&
&
\coho H * {M}  \ar[ur]  |\2 \ar[dr] |\1 &
&
\coho H{ *+1} {M/\stres,{M^{\stres}} }  \ar[ur] |\4\ar[dr] |\3 &
&
\coho H {*+1} {M}  \ar[ur] |\1\ar[dr] |\2&
\\
&
\coho H * {M,{M^{\stres}} }  \ar[ur] |\2\ar@/_{7mm}/[rr]  |\4&
&
\coho H {*-3} {M/S^3, \Sigma/S^3}  \oplus \left( \coho H {*-2} {M^{\sbat}}\right)^{-\zdos} \ar[ur] |\4\ar@/_{7mm}/[rr] |\1&
&
\coho H {*+1} {M/\stres}  \ar[ur] |\1\ar@/_{7mm}/[rr] |\3&
&
\coho H {*+1} {{M^{\stres}} } .
}}
\eee

\smallskip

Furthermore, using property \eqref{B1}, we obtain the long exact sequences
$$
\xymatrix@C=0,35cm
{
\cdots \ar[r] &
\coho H {*} M 
\ar[r]  &\coho H {*-3} {M/S^3, \Sigma/S^3}  \oplus \left( \coho H {*-2} {M^{\sbat}}\right)^{-\zdos} \oplus  \coho H {*} {M^{S^3}} 
\ar[r]  &
 \coho H {*+1} {M/S^3,M^{S^3}} 
\ar[r]  &
 \coho H {*+1} M
\ar[r]  &\cdots,}
$$
the Smith-Gysin sequence, and
$$
\xymatrix@C=0,35cm
{
\cdots \ar[r] &
 \coho H {*}{ M/\stres, M^{\stres}}  \ar[r]  &   
\coho H {*} {M/{S^3}}  \oplus  \coho H {*} { M, M^{\stres}}  \ar[r]  &
 \coho H {*} {M} \ar[r]  &
 \coho H {*+1}{ M/\stres, M^{\stres}} \ar[r]  & 
 \cdots,
}
$$
We have constructed the two long exact sequences.
\epr

\begin{remarque}

\nada 

1 - The technique used to construct the Smith-Gysin sequence is the same as the one used in \cite[page 161]{MR0413144}, where a smooth action\footnote{The result obtained by Bredon is more general than the one presented in this Remark, which only applies to smooth manifolds.} $\Phi \colon S^1 \times M \to M$ is given, and the following Smith-Gysin sequence is constructed:
\SucExa
{
\coho H {*} M
}{
 \coho H {*-1} {M/S^1, M^{S^1}} \oplus \coho H {*} {M^{S^1}}
}{
\coho H {*+1} {M/S^1, M^{S^1}}
}{
 \coho H {*+1} M }
It is worth noting that the long exact sequence
\SucExaRef{S1SmGy}
{
 \coho H {*}{ M/\sbat, M^{\sbat}}
}{
\coho H {*} {M/{S^1}}  \oplus  \coho H {*} { M, M^{\sbat}}
}{ \coho H {*} M 
}{
 \coho H {*+1}{ M/\sbat, M^{\sbat}} 
 }
can be also be obtained. 

Recall that the associated Gysin sequence is 
\SucExaRef{S1Gys}
{\coho H {*} M} {\coho H {*-1} {M/S^1, M^{S^1}} } {\coho H {*+1} {M/S^1} } {\coho H {*+1} M }
(cf. \cite[page 161]{MR0413144}).
\smallskip

2-  We can obtain the Gysin sequences  \eqref{gys3SemFr}, \eqref{S1Gys} from the Smith-Gysin sequences \eqref{gys3ClasSG},\eqref{S1SmGy} through the following steps.
		
Assume that $G=S^k$ acts semi-freely on $M$, where $k$ is either 1 or 3, and let $F$ be the set of fixed points. The quotient map $\pi\colon M\longrightarrow B=M/S^k$ is continuous, and its restriction to $M\setminus F$ is a $G$-principal bundle over $B\setminus F$, where we have identified $F$ with $\pi(F)$ for convenience. This restriction induces the pullback map 
$\pi^*\colon \Omega_v^*(B)\to\Omega_v^*(M)$ iin the Verona complexes, which, in turn, induces a map in cohomology that we denote by $\pi^*\colon H^*(B)\to H^*(M)$, using isomorphism in cohomology.
 We have the following diagram
$$
\scalebox{.9}
{
\xymatrix@C=1.8cm{
	0\ar[d] & 0\ar[d] & 0\ar[d] & 0\ar[d]  \\
	H^i(B,F)\ar[r]^{\iota}\ar[d]^{(1,\iota)} &
		 H^i(B)\ar[r]^{\rho_B}\ar[d]^{(\pi^*,1_B)} & 
			 H^i (F)\ar[r]^{\delta}\ar[d]^{(0,1_F)} &
				H^{i+1}(B,F)\ar[d]^{(1,\iota)} \\
	H^i(B,F)\oplus H^i(B)\ar[r]^{(\pi^*\circ\iota,1_B)}\ar[d]^{\iota-1_B} &
		 H^i(M)\oplus H^i(B)\ar[r]^{(\int,\pi^*\circ\rho_M)}\ar[d]^{1_M-\pi^*} & 
			 H^{i-k}(B,F)\oplus H^i (F)\ar[r]^{(\wedge[e]+\delta,0)}\ar[d]^{proj_1} & 
				 H^{i+1}(B,F)\oplus H^{i+1}(B)\ar[d]^{\iota-1_B} \\
	 H^i(B)\ar[r]^{\pi*}\ar[d] & 
		 H^i(M)\ar[r]^{\int}\ar[d] & 
			 H^{i-k}(B,F)\ar[r]^{\wedge[e]}\ar[d] & 
				 H^{i+1}(B)\ar[d]\\
	 0 & 0 & 0 & 0 }
	 }
$$
The top horizontal sequence is the long exact sequence of the pair $(B,F)$. The middle sequence is the direct sum of the Smith-Gysin sequence of \cite{Smith-Gysin} and the exact sequence 
$$
\cdots\to 0\to H^i(B)\stackrel{1_B}{\longrightarrow} H^i(B)\to 0 \to H^{i+1}(B)\stackrel{1_B}{\longrightarrow} H^{i+1}(B)\to0\to \cdots,
$$
which is exact as well. The above diagram is commutative, and we can verify its commutativity easily. By using the commutativity of the diagram, the exactness of the two horizontal sequences, and the vertical short sequences, we can perform a diagram chase that shows the exactness of the bottom long sequence (known as the Gysin sequence).
\end{remarque}

\end{document}